\documentclass{amsart}
\usepackage{amsmath,amsthm,amsfonts,amssymb,amscd,mathrsfs}
\usepackage{psfrag,graphicx}

\usepackage{epic,eepic}
\usepackage{color}



\theoremstyle{plain}

\newtheorem{Thm}{Theorem}[section]
\newtheorem{Lem}[Thm]{Lemma}

\newtheorem{Cor}[Thm]{Corollary}
\theoremstyle{remark}

\long\def\begcom#1\endcom{}

\newcommand{\length}{\operatorname{\length}}

\def\length{\operatorname{length}}

\email{chaol@cufe.edu.cn}

\begin{document}
\begin{center}
\title[On Density Of Positive Lyapunov Exponents For $C^1$ Symplectic Diffeomorphisms]
{On Density Of Positive Lyapunov Exponents For $C^1$ Symplectic Diffeomorphisms}

\thanks{$^{*}$ CL is supported by NSFC(\# 11471344), Beijing Higher Education Young Elite Teacher Project and Program for Innovation Research in Central University of Finance and Economics}
\medskip

\date{April, 2014}

\maketitle

\centerline{\scshape Chao Liang}
\medskip
{\footnotesize
 \centerline{Applied Mathematical Department}
   \centerline{Central University of Finance
and Economics}
   \centerline{ Beijing, 100081, China}
}

\begin{abstract}
Let $M$ be a 2$d-$dimensional compact connected Riemannian manifold and $\omega$ be a symplectic form on $M$. In this paper, we prove that a symplectic diffeomorphism, with all Lyapunov exponent zero for almost everywhere, can be $C^1$ approximated by one with a positive Lyapunov exponent for a positive-measured subset of $M$. That is,
the set
\[
\left\{
f\in \mathcal{S}ym^1_{\omega}(M)\,|
\begin{array}{ll}
&\mbox{The largest Lyapunov exponent }\lambda_1(f,\,x)>0\\
&\mbox{ for a positive measure set }
\end{array}
\right\}
\]
is dense in
$\mathcal{S}ym^1_{\omega}(M)$.
\end{abstract}
\end{center}
\bigskip

\section{INTRODUCTION}

Lyapunov exponent is a useful tool to describe hyperbolicity of  a system. Using this concept, Pesin introduced a weaker form of hyperbolicity, which he termed nonuniform hyperbolicity. A system that admits an invariant measure with nonzero Lyapunov exponents for almost everywhere is nonuniformly hyperbolic. Since uniformly hyperbolic systems are not dense in the set of smooth dynamical systems, Pesin questioned that whether non-uniformly hyperbolic systems are dense or generic in $C^r$ volume-preserving diffeomorphisms. KAM theory gives a negative answer to this conjecture when $r$ is large enough. And later Ma$\tilde{n}\acute{e}$-Bochi-Viana's results\cite{B,BV} imply that nonuniformly hyperbolic systems can not be generic. Hence Pesin's question should be whether non-uniformly hyperbolic systems are dense in $C^1$ volume-preserving diffeomorphisms. Thus it is of utmost importance to detect when the zero exponents can be removed by perturbations.
Shub-Wilkinson's example\cite{SW} builds a conservative
perturbation to a skew product of an Anosov diffeomorphism of the
torus $\mathcal{T}^2$ by rotations and create positive exponents in
the center direction for Lebesgue almost every point.
Baraviera-Bonatti\cite{BB} present a local version of Shub-Wilkinson's
argument, allowing one to perturb the sum of all the center
integrated Lyapunov exponents of any conservative partially
hyperbolic systems. These perturbations can be assumed to be done on some invariant bundles (not all bundles) of a dominated splitting by decreasing the larger integrated Lyapunov exponent and at the same time increasing the smaller one. But the changes should be different, which is exactly the new integrated Lyapunov exponent of the center bundle. For symplectic difeomorphisms, one can not do this kind of perturbation since the absolute values of the changes will always be the same. Thus, we want to consider the question that for a system with all the Lyapunov exponents zero, which implies that we can not `borrow' some nonzero exponent from one bundle to perturb zero exponent on the center bundle, wether it can be approximated by a new system with positive Lyapunov exponents. That is the density of systems with positive Lyapunov exponents for a positive-measured set. This paper is to deal with the problem.

Let $M$ be a 2$d-$dimensional compact connected Riemannian manifold and $\omega$ be a symplectic form on $M$, i.e. a non-degenerate closed $2-$form. Taking $d$
times the wedge product of $\omega$ with itself we obtain a volume form on $M$. A $C^r$
diffeomorphism $f$ of $M$, $r\geq 1$, is called symplectic if it preserves the symplectic
form, $f^*\omega=\omega$. Denote by $\mathcal{S}ym^1_{\omega}(M)$ the set of all $C^1$ symplectic diffeomorphisms on $M$. Our main result is the following theorem.
\smallskip

\begin{Thm}\label{ThmPosLyaExp}
The set $\{f\in \mathcal{S}ym^1_{\omega}(M)\,|\,\mbox{The largest Lyapunov exponent }\lambda_1(f,\,x)>0
\mbox{ for a positive measure set }\}$ is dense in
$\mathcal{S}ym^1_{\omega}(M)$.
\end{Thm}
\smallskip

\section{PROOF OF THEOREM \ref{ThmPosLyaExp}}

In 1979, Katok\cite{K} gave an example of 2-dimensional nonuniformly hyperbolic system. It is derived from an Anosov map by a big $C^1-$perturbation. We use the product of finitely many Katok's examples to construct a nonuniformly hyperbolic system of high dimension. In this section, we will prove our main result by pasting this example to neighborhoods of elliptic periodic points. First, let us introduce some notations and an important lemma as the following.

If for a periodic point $p$ of period $k$ the tangent map $Df^k(p)$ has exactly $2l$ simple
non-real eigenvalues of norm $1$ and the other ones have norm different from $1$,
then we say that $p$ is {\em an $l-$elliptic periodic point}. Sometimes it is also called quasi-elliptic.

\begin{Lem}\label{Lem1}
  Let $f\in \mathcal{S}ym^1_{\omega}(M)$ with all Lyapunov exponents zero. And $p$ be an elliptic periodic point of $f$. Then there is an arbitrarily small perturbation $g\in \mathcal{S}ym^1_{\omega}(M)$ of $f$ with a positive Lyapunov exponent on a positive-measured set.
\end{Lem}


{\bf Proof}\,Let $n$ be the period of the periodic point $p$. By a small
perturbation, we can assume that $D^nf$ is a rational rotation. Then
there exists a positive integer $N>0$ and a measurable set $B$ in
$M$, such that $f^N|_B=id$.

In the following we will construct a perturbation $\hat{h}$ of $id$
and define a perturbation $g$ of $f$ such that $g=f$ outside
$\bigcup_{1\leq i\leq N} f^iB$ and $g=\hat{h}\circ f$ in
$\bigcup\limits_{\ell}\bigcup\limits_{k} R_{2\pi j/k}D^\ell$.

Take $h=\underbrace{K\times\cdots\times K}_{d \mbox{ times}}$, where $K$ is Katok's example in \cite{K}. Then $h$ is volume preserving but not a small
perturbation of $id$. Find a flow $\phi_t$, an isotopic to
$\mathbb{D}$ satisfying that
$$\phi_t:\,B\times[0,\,1]\rightarrow B$$
$$\phi_0=id,\quad \phi_1=h.$$
Take $k$ large and define a sequence of small perturbation of
$id:\,\{h_i\}^{k}_{i=1}$ such that
$$h_i\circ\cdot\cdot\cdot\circ h_0=\phi_{i/k},\,\forall i=0,1,...,k-1.$$
Then $h_i$ is volume preserving and $h=h_k\circ\cdot\cdot\cdot\circ
h_0$. Since $k$ is large, the rotation $R_{2\pi/k}$ is a small
rotation with the angle $2\pi/k$. Find a small disk $D^\ell\subset
f^\ell B$, such that
$$R_{2\pi/k}D^\ell\bigcap D^\ell=\emptyset,\,\forall \ell=0,1,...,N-1.$$
Define
\[
\hat{h}(x)=\left\{
\begin{array}{ll}
h_{j-1}\circ R_{2\pi/k}\circ\cdot\cdot\cdot\circ R_{2\pi/k}\circ
h_1\circ R_{2\pi/k}&\quad \mbox{ on } R_{2\pi j/k}D^\ell\\
id &\quad \mbox{ outside } f^\ell B
\end{array}
\right.
\]
$i=1,...,k.$
\bigskip

A $Df-$invariant splitting $TM=E^1\oplus\cdots\oplus E^k$ is called
a {\em dominated splitting} if each $E^i$ is a continuous
$Df-$invariant subbundle of $TM$ and if there is some integer $n>0$
such that, for any $x\in M$, any $i<j$ and any non-zero vectors
$u\in E^i(x)$ and $\nu\in E^j(x)$, one has
$$\frac{\|Df^n(u)\|}{\|u\|}<\frac{1}{2}\frac{\|Df^n(\nu)\|}{\|\nu\|}.$$
A map $f$ is called {\em partially hyperbolic} if there exists a dominated splitting $TM
= E^u\oplus E^c\oplus E^s$, into nonzero bundles such that,
for some Riemannian metric $\|\cdot\|$ on M, we have $$\|(Df
|E^u(x))^{-1}\|^{-1} > \|Df |E^c(x)\|\geq\|(Df |E^c(x))^{-1}\|^{-1}
> \|Df |E^s(x)\|$$ for every $x\in M$. Such a splitting is automatically
continuous.
where $E^s$ and $E^u$ denote the strong expanding and contracting
invariant bundle, respectively.

If a system if partially hyperbolic, then at least one of the Lyapunov exponents is nonzero. So we consider the systems that are not partially hyperbolic. The existence of an elliptic
periodic point is an obstruction for partial hyperbolicity. Saghin-Xia\cite{SX} proved that the
converse is also true generically, i.e. if a $C^1-$generic symplectic diffeomorphism is
not partially hyperbolic, then it has an elliptic periodic point (actually it has a
dense set of elliptic periodic points).

\begin{Lem}[Theorem 1 in \cite{SX}]\label{LemSX}
There exists an open dense subset $\mathcal{U}$ of $\mathcal{S}ym^1_{\omega}(M)$ such that any function
in $\mathcal{U}$ is either partially hyperbolic or it has an elliptic periodic point. There
exists a residual subset $\mathcal{R}$ of $\mathcal{S}ym^1_{\omega}(M)$ such that any function in $\mathcal{R}$ is either partially
hyperbolic or the set of elliptic periodic points is dense on the manifold.
\end{Lem}
\medskip

\begin{proof}[\bf Proof of Theorem \ref{ThmPosLyaExp}]\,We will prove by induction. For every $n\in \mathbb{N}$, by Lemma \ref{LemSX}, we can find finitely many elliptic periodic points such that they form an $\frac1n-$net of $M$. We can find small neighbourhoods of them such that they are disjoint from each other and also disjoint from all neighbourhoods we find in the preceding steps for integers smaller than $n$. And then we can do the perturbation as Lemma\ref{Lem1} in these neighbourhoods. Then we get our result.
\end{proof}
\smallskip

By the proof above, we can further obtain the following corollary.

\begin{Cor}
Let $M$ be a compact surface without boundary. The set
\[
\left\{
f\in \mathcal{S}ym^1_{\omega}(M)\,|
\begin{array}{ll}
&\mbox{The largest Lyapunov exponent }\lambda_1(f,\,x)>0\\
&\mbox{ for an open, dense and positive-measured set }
\end{array}
\right\}
\]
is dense in
$\mathcal{S}ym^1_{\omega}(M)$.
\end{Cor}
\medskip

Let $M$ be a compact surface without boundary. The symplectic form $\omega$ is also a volume form. By Oseledets Theorem, on almost all points, the Lyapunov exponents exist and the sum of the two exponents are zero. We say that $f$ is {\em nonuniformly hyperbolic} on a point if the two Lyapunov exponent are nonzero along the orbit of this point.

\begin{Cor}
Let $M$ be a compact surface without boundary. The set
\[
\left\{
f\in \mathcal{S}ym^1_{\omega}(M)\,|
\begin{array}{ll}
&f\mbox{ is nonuniformly hyperbolic }\\
&\mbox{ on an open, dense and positive-measured set }
\end{array}
\right\}
\]
is dense in
$\mathcal{S}ym^1_{\omega}(M)$.
\end{Cor}

{\bf Acknowledgement}\,\,We want to thank Professor Bochi and Professor Zhihong Xia.

\bigskip

\end{document}